\documentclass[10pt, reqno]{amsart}

\newcommand{\PD}{\partial}

\newcommand{\wt}{\tilde} % modified

\DeclareMathOperator{\WF}{WF_A}

\newcommand{\Ac}{\mathcal{A}}

\newcommand{\Cc}{\mathcal{C}}
\newcommand{\Dc}{\mathcal{D}}
\newcommand{\Ec}{\mathcal{E}}

\newcommand{\Uc}{\mathcal{U}}

\newcommand{\Rb}{\mathbb{R}}
\newcommand{\Sb}{\mathbb{S}}

\newcommand{\R}{\rangle}
\newcommand{\Beq}{\begin{equation}}
\newcommand{\Eeq}{\end{equation}}
\newcommand{\beq}{\begin{equation*}}
\newcommand{\eeq}{\end{equation*}}
\newcommand{\bal}{\begin{align}}
\newcommand{\eal}{\end{align}}

\renewcommand{\L}{\langle}

\newcommand{\g}{\gamma}

\newcommand{\A}{\alpha}

\newcommand{\bp}{\begin{prob}}
\newcommand{\ep}{\end{prob}}
\newcommand{\bpr}{\begin{proof}}
\newcommand{\epr}{\end{proof}}

\newcommand{\inter}{^\text{\rm int}}

\usepackage{amssymb,amsmath,amsfonts,amsthm}
\newtheorem{thm}{Theorem}
\newtheorem{lemma}{Lemma}
\newtheorem{prop}{Proposition}
\newtheorem{defn}{Definition}

\theoremstyle{definition}

\newtheorem{prob}{Problem}

\newtheorem{remark}{Remark}
\usepackage{amsmath}
\begin{document}
\title[A support theorem for the geodesic ray transform]{A support theorem for the geodesic ray transform of symmetric tensor fields}

\author[V. Krishnan]{Venky Krishnan}
\address{Department of Mathematics, Tufts University, Medford, MA 02155}

\author[P. Stefanov]{Plamen Stefanov}
\address{Department of Mathematics, Purdue University, West Lafayette, IN 47907}
\thanks{Second author partly supported by NSF}

\begin{abstract}
Let $(M,g)$ be a simple Riemannian manifold with boundary and consider the geodesic ray transform of symmetric 2-tensor fields. Let the integral of $f$  along maximal geodesics vanish on an appropriate open subset of the space of geodesics in $M$. Under the assumption that the metric $g$ is real-analytic, it is shown that there exists a vector field $v$ satisfying $f=dv$ on the set of points lying on these geodesics and $v=0$ on the intersection of this set with the boundary $\PD M$ of the manifold $M$. Using this result, a Helgason's type of a support theorem for the geodesic ray transform is proven. The approach is based on  analytic microlocal techniques.
\end{abstract}

\maketitle
\section{Introduction}\label{S1}
Let $(M,g)$ be an $n$-dimensional smooth compact Riemannian manifold with boundary $\PD M$. We define the geodesic ray transform of a symmetric $2$-tensor field $f$ %\in L^{2}(M)$ 
as
\Beq
If(\g)=\int\limits_{0}^{l(\g)}f_{ij}(\g(t))\dot{\g}^{i}(t)\dot{\g}^{j}(t)dt,
\Eeq
where $\g:[0,l(\g)]\to M$ is any geodesic parameterized by its arc length with  endpoints on the boundary $\PD M$.

The transform $I$ is not injective. Indeed, any potential tensor $dv$ with $v=0$ on $\PD M$ belongs to its kernel \cite{Sh1}. Here $dv$ is the symmetrized covariant derivative of  the one-form $v$ given locally by $(dv)_{ij} = \frac12(\nabla_iv_j+\nabla_jv_i)$. It is expected that this is the only obstruction to injectivity for certain classes of manifolds, including simple ones. We call this s-injectivity of $I$. Simple manifolds are those that are  convex,  have no conjugate points and have strictly convex boundary, see next section.  

The geodesic ray transform $I$ arises naturally in boundary and lens rigidity problems, see e.g., \cite{Sh1, SU1, SU2, SU-AJM} and the references there. Recall that $(M,g)$ is called boundary rigid, if it is determined uniquely, up to isometry, by the boundary distance function $\rho_g(x,y)$, $(x,y)\in \PD M\times \PD M$. A linearization of $\rho$ near a fixed simple metric (and $M$ fixed) is given by $If$, where $f$ is the variation of $g$. So the problem of inverting $I$ modulo potential fields is a linear version of the problem of recovering $g$ from $\rho_g(x,y)$, $(x,y)\in \PD M\times \PD M$ modulo isometries. 

In this paper, we prove support type of theorems for $I$ on simple real-analytic  Riemannian manifolds.  Without any conjugacy assumptions, $I$ may fail to be s-injective, as a simple example based on the symmetry of the sphere shows \cite{SU-AJM}. On the other hand,  there are  non-simple manifolds, for which we still have \mbox{s-injectivity}, see \cite{SU-AJM}.

\section{Main Results} 
\begin{defn}
We say that a compact Riemannian manifold $(M,g)$ with boundary is simple if
\begin{enumerate}
\item[(a)] The boundary $\PD M$ is strictly convex: $\L \nabla_{\xi}\nu,\xi\R > 0$ for each $\xi \in T_{x}(\PD M)$ where $\nu$ is the unit outward normal to the boundary.
\item[(b)] The map $\exp_{x}:\exp_{x}^{-1}M\to M$ is a diffeomorphism for each $x\in M$.
\end{enumerate}
\end{defn}

Note that (b) also implies that any two points in $M$ are connected by unique geodesic depending smoothly on its endpoints. This, together with (a) and the assumption of no conjugate points is an equivalent characterization of a simple manifold. 
Simple manifolds are diffeomorphic to balls, see \cite{Sh-UW} 

From now on, we assume that $M$ is a simple  manifold. We fix a real analytic atlas on $M$. Actually, just one chart would be enough, so we can think of $M$ as the closure of a  subdomain of $\Rb^n$. 
We  say that a function, or more generally, a tensor field is analytic on the set $U$, not necessarily open, if it is real analytic in some neighborhood of $U$. Let $g$ be an  analytic Riemannian metric  on $M$. In particular, we have that $M$  can be extended  to a slightly larger  simple real analytic manifold $\wt{M}$ such that %$\wt{M}$ is simple and 
$M$ is in the interior of $\wt{M}$, and  $g$ extends  to $\wt{M}$ to a  metric (which we still denote by $g$) that is still real analytic. Let us also extend all symmetric tensor fields $f$, a priori defined on $M$ only, as $0$ in $\wt{M}\setminus M$.

We can parametrize maximal geodesic segments (that we call geodesics) by their endpoints on the boundary $x$, $y$. We will use the notation $\gamma_{[x,y]}$ for that. Then $\gamma_{[x,x]}$ is just a point. This parametrization induces a natural topology in the set of all geodesics.

Given a set $\Ac$ of geodesic, we denote by $M_{\Ac}$ the set of points lying on the geodesics in $\Ac$. Also we denote by $\PD_{\Ac}M$  the intersection of $M_{\Ac}$ with the boundary $\PD M$. We extend all geodesics to $\tilde M$, and we call the set $\Ac$ again. 

We extend $I$  by duality to act on tensor fields that are distributions in $\tilde M$ supported in $M$.  We denote the space of compactly supported tensors by $\Ec'(\tilde M)$, and it will be clear from the context what the order of the tensor is. The condition $v=0$ on $\PD M$ then is replaced by the condition that the distribution-valued 1-form $v$ vanishes outside $M$. This is consistent with the classical case. Indeed, if $f\in C^1(M)$, and if $f=dv$ with $v=0$ on $\PD M$, then $f=dv$ remains true for $f$ and $v$ extended as zero to $\tilde M\setminus M$ because the operation of extension as zero and the differential $d$ commute on $v$'s vanishing on $\PD M$. 

We now give the statements of our main results.

\begin{thm}\label{S1:T1}
Let $(M,g)$ be a simple analytic Riemannian manifold. 
Assume that  $\Ac$ is an open set of geodesics satisfying
\begin{equation}   \label{i}
\mbox{any geodesic in $\Ac$ is homotopic, within the set $\Ac$, to a point on $\PD M$.}
\end{equation}

(a) Then, given a symmetric $2$-tensor field $f\in \Ec'(\tilde M)$ supported in $M$, we have that $If(\g)=0$ for each $\g\in \Ac$ if and only if  for each $\g_0\in \Ac$, 
there exists a neighborhood $\Uc$ of $\g_0$ and an 1-form $v\in \Dc'(\tilde M_{\Uc})$,  so that $f=dv$ in $\tilde M_{\Uc}$, and $v=0$ outside $M$.

(b) If, in addition, 
\begin{equation}   \label{ii}
\mbox{$\pi_{1}(M_{\Ac},\PD_{\Ac} M)=0$,}
\end{equation}
then there is a unique $v$ defined in $\tilde M_\Ac$, vanishing outside $M$ so that $f=dv$ in $\tilde M_\Ac$. 
\end{thm}

Condition \eqref{i} means the following: Given $\gamma=\gamma_{[x,y]}\in \Ac$, there exists continuous curves $\alpha(t)$, $\beta(t)$, $0\le t\le1$, on $\PD M$, so that (i) $\alpha(0)=x$, $\beta(0)=y$; (ii) $\gamma_{[\alpha(t),\beta(t)]}\in \Ac$; and (iii) $\alpha(1)=\beta(1)=:z$ and  $\gamma_{[z,z]}\in \PD M$. 

Condition \eqref{ii} means that any closed path with a base point on $\PD M$ is homotopic to a path lying on $\PD M$, see also \cite{AH, LeeU}.

Our next theorem is a version of Helgason's support theorem in the geodesic ray transform setting.

\begin{thm}\label{S1:T2}
Let $(M,g)$ be a simple analytic Riemannian manifold. 
 Let $K$ be a closed geodesically convex subset of $M$. If for a symmetric $2$-tensor field $f\in \Ec'(\tilde M)$ supported in $M$, we have that 
$If(\g)=0$ for each geodesic $\g$ not intersecting $K$, then there exists an 1-form $v\in \Dc'(\tilde M^{\inter}\setminus K)$  such that $f=dv$ in $ \tilde M^{\inter}\setminus K $, and $v=0$ in $\tilde M^{\inter}\setminus M$.
\end{thm}

We say that $K$ is geodesically convex if for any two points in $K$, the unique  geodesic that connects them lies in $K$ as well. Here and below, $B^{\inter}$ stands for the interior of the set $B$. 

Although we have stated our results only for symmetric $2$-tensor fields, these results hold with minor modifications for tensor fields of any rank. For ease of notation and readability, we have limited ourselves to tensors of rank $2$ only.  Similar results for functions (tensors of rank $0$) were obtained by the first author in \cite{Kr}.

It is well known than a symmetric 2-tensor field $f\in L^{2}(M)$ can be written uniquely as the orthogonal sum of two fields; a solenoidal part $f^{s}$ and a potential part $dv$ with $v|_{\PD M}=0$ \cite{Sh1}. 
Then $I(dv)=0$ is a consequence of the following identity
\begin{equation}\label{dv}
\frac{d}{dt}\langle v(\gamma(t)), \dot\gamma(t)\rangle = 
\langle dv(\gamma(t)) , \dot\gamma(t)\otimes \dot\gamma(t)\rangle.
\end{equation}
Following \cite{SU1}, we say $I$ is $s$-injective, if $If=0$ implies $f^{s}=0$. S-injectivity under a small curvature assumption has been established in \cite{Sh1, Sh-UW, Dairbekov}. In the 2D case, s-injectivity is known for all simple metrics \cite{Sh-2D}. 
The second author and Uhlmann have studied the question of s-injectivity of geodesic ray transform in \cite{SU1, SU2, SU-AJM}. Among the results there, it is shown that on a compact Riemannian manifold with boundary, for a generic set of simple metrics that includes real-analytic simple metrics, the geodesic ray transform is $s$-injective. Also in \cite{SU-AJM}, it is  shown that the same result is true for a class of non-simple Riemannian manifolds and for the case where the geodesic ray transform is restricted to certain subsets $\Gamma$ of the space of geodesics. The set $\Gamma$ is such that the collection of their conormal bundles covers $T^{*}M$.

This paper also deals with the partial data case as in \cite{SU-AJM}, but we allow for a more general open set $\Gamma$. In particular, we do not impose the condition that the collection of conormals of $\Gamma$ covers $T^{*}M$. Thereby, we are able to prove support theorems for the geodesic ray transform.

We use analytic microlocal analysis to prove our results. Guillemin first introduced the microlocal approach in the Radon transform setting. For more details see the book \cite{GS}. In 1987, Boman and Quinto in \cite{BQ1} used analytic microlocal techniques to prove support theorems for Radon transforms with positive real-analytic weights. Since then, other support theorems have been proven by Boman, Grinberg, Gonzalez, Quinto and Zhou in both collaborative and individual works. Some references are \cite{B1, B2, BQ1, BQ2, GQ1, GQ2, Q1, Q2, QZ}. Also the results of the second author and Uhlmann in \cite{SU2, SU-AJM} involve  analytic microlocal analysis. Our proofs rely on techniques introduced in 
\cite{SU1, SU2, SU-AJM}.

\section{Preliminary constructions}\label{S2} As we mentioned above, we extend all geodesic in $\Ac$ to maximal geodesics on $\tilde M$ to obtain a new set denoted by  $\Ac$ again. We also add to $\Ac$ all geodesics in $\tilde M$ that do not intersect $M$. Clearly, the new set $\Ac$ is open again, and \eqref{i} is preserved, with $M$ replaced by $\tilde M$.

Let $\gamma_0$,  be a maximal geodesic in $\tilde M$ connecting $x_0\not= y_0$. Using normal coordinates at $x_0$, one can easily construct coordinates $x=(x',x^n)$ in $\tilde M\setminus \{x_0\}$ so that $x^n$ is the distance  to $x_0$, and $\partial/\partial x^n$ is normal to $\partial/\partial x^\alpha$, $\alpha<n$, see \cite{SU1}. They are actually semi-geodesic normal coordinates to any geodesic sphere centered at $x_0$. 
In those coordinates, the metric $g$ satisfies $g_{ni}=\delta_{ni}$, $\forall i$, and hence the Christoffel symbols satisfy $\Gamma_{nn}^{i}=\Gamma_{ni}^{n}=0$. 
Thus we get global coordinates in $\tilde M\setminus \{x_0\}$, and in particular, in $M$, of the type above. All lines of the type $x'=\mbox{const.}$ are now geodesics with $x^n$ an arc length parameter. 

Let $U$ be a tubular neighborhood of $\gamma_0$ in $M$, given by $|x'|<\epsilon$, $a(x')\le x^n\le b(x')$, where $a$ and $b$ are analytic functions giving a local parametrization of $\partial M$.  
Given a symmetric 2-tensor field $f$ (for now, we assume $f\in C^\infty(M)$), one can always construct an one form $v$ in $U$ so that for 
\begin{equation}   \label{ha}
h :=f-dv
\end{equation}
one has
\begin{equation}  \label{h}
h_{ni}=0,\quad  \forall i,\quad v(x',a(x'))=0.
\end{equation}
The form $v$ is defined by solving a system of  ODE's. This construction can be found in several papers, see \cite{Eskin, Sh-sibir, Sh1, SU1, SU2}, and is based on the following.  We have $\nabla_{i}v_{j}=\PD_{i}v_{j}-\Gamma_{ij}^{k}v_{k}$. First let $i=n$. Then we solve 
\[ (dv)_{nn}=\PD_{n}v_{n}-\Gamma_{nn}^{k}v_{k}=f_{nn}
\]
for $v_n$. 
Since in $U$, $\Gamma_{nn}^{k}=0$, we get
\begin{equation}\label{h2}
\PD_{n} v_{n}=f_{nn},\quad v_{n}(x',a(x'))=0.
\end{equation}
Using this we solve for the remaining $v_{i}$'s. We have $\PD_{n}v_{i}+\PD_{i}v_{n}-2\Gamma_{ni}^{k}v_{k}=2f_{in}$. Since $\Gamma_{ni}^{n}=0$ for $0\leq i\leq n$, this set of equations reduces to a linear system for the variables $v_{i},1\leq i\leq i-1$:
\begin{equation}  \label{h3}
\PD_{n}v_{i}-2\Gamma_{ni}^{\A}v_{\A}=2f_{in}-\PD_{i}v_{n},\quad\!\! v_{i}(x',a(x'))=0,\quad\!\! 1\leq \A\leq n-1, \quad\!\! 1\leq i\leq n-1.
\end{equation}
If we assume that $If=0$, then we get that $v_n$ vanishes on the other part of $\partial M$, as well, i.e., for $x^n=b(x')$. We do not automatically get this for $v^i$, $i<n$ though. Next lemma shows that this is actually true for those components of $v$ as well.

We denote now by $\tilde U$ a tubular neighborhood of $\gamma_0$ of the same type but related to $\tilde M$.  
\begin{lemma} \label{lemma_v0}
Let $f$ be supported in $M$, and $If(\gamma)=0$ for all maximal geodesics in $ \tilde U$ belonging to some neighborhood of the geodesics $x'=\mbox{\rm const}$.  Then  $v=0$  in $ \tilde U\inter\setminus M$. 
\end{lemma}

\bpr
Let first $f\in C^\infty(M)$. 
We will give another, invariant definition of $v$. For any $x\in \tilde U$ and any $\xi\in T_x \tilde U\setminus \{0\}$ so that $\gamma_{x,\xi}$ stays in $ \tilde U$, we set
\begin{equation} \label{u}
u(x,\xi) = \int_{\tau_-(x,\xi)}^0 f_{ij}(\gamma_{x,\xi}(t))\dot \gamma_{x,\xi}^i(t) \dot\gamma_{x,\xi}^j(t)\, dt,
\end{equation}
where $\tau_-(x,\xi)\le0$ is determined by $\gamma_{x,\xi}( \tau_-(x,\xi)  )\in \PD M$. Extend the definition of $\gamma_{x,\xi}$ for $\xi\not=0$ not necessarily unit as a solution of the geodesic equation. 
Then $u(x,\xi)$ is homogeneous of order $1$ in $\xi$. Moreover, $u$ is odd in $\xi$ for $\xi/|\xi|$ close enough to $\pm e_n :=\pm (0,\dots,0,1)$ because $If=0$ there.  If $f=dv$ with $v=0$ outside $M$, then $u(x,\xi)=\langle v(x), \xi\rangle$ with $v$ considered as a vector field here. This is the basis for our definition.  We set (now $v$ is identified with a covector field using the metric)
\begin{equation}  \label{vv}
v(x) = \partial_\xi\big|_{\xi=e_n} u(x,\xi), \quad \forall x\in \tilde U.
\end{equation}
It is easy to check that $v=0$  outside $M$, i.e., for $x^n<a(x')$ or for $x^n>b(x')$. In particular, we get that \eqref{h2} is satisfied because $\xi^k \partial u/\partial \xi^k = u$ by the homogeneity of $u$; thus setting $\xi=e_n$, we get $v_n(x) = u(x,e_n)$.  
We need to show that for $h=f-dv$, one has $h_{ni}=0$, see  \eqref{h}. Let $u$ be as above but related to $h$ instead of $f$. Then
\begin{equation}\label{h4}
u(x,e_n)=0, \quad \partial_\xi\big|_{\xi=e_n} u(x,\xi)=0.
\end{equation}
We also have $Gu(x,\xi)=h_{ij}\xi^i\xi^j$, where $G=\xi^i\partial_{x^i}-\Gamma^k_{ij}\xi^i\xi^j\partial_{\xi^k}$ is the generator of the geodesic flow. Differentiate w.r.t.\ $\xi^\alpha$, $\alpha<n$, at $\xi=e_n$, using \eqref{h4} and the fact that $\Gamma_{nn}^k=0$ to get $0=h_{\alpha n}$. Since the field $v$ with the properties \eqref{h} is unique, this completes the proof of the lemma in the case $f\in C^\infty$.

Let now $f$ be a distribution in $\tilde M$, supported in $M$.  We claim that $u(x,\xi)$ is a well defined distribution of $x$ (in the interior of $\tilde M$) depending smoothly on $\xi$. Indeed, let $\phi(x)$ be a test function supported in the interior of $\tilde M$. We can always assume that its support is small enough. The map $x\mapsto y=\gamma_{x,\xi}(t)$ is a local diffeomorphism depending smoothly on $t$ and $\xi$ for $\xi$ close enough to $e_n$. Indeed, this is true for $\xi=e_n$ because then $\gamma_{x,\xi}(t)=x+t\xi$; and for $\xi$ close to $e_n$ it is true by a perturbation argument. Let $x=x(y,\xi,t)$ be the resulting function. Then we can also write $\eta(y,\xi,t) := \dot \gamma_{x,\xi}^i(t)$. 
Then
\[
\begin{split}
\int u(x,\xi)\phi(x)\, dx &= \int\int_{-\infty}^0  f_{ij}(\gamma_{x,\xi}(t))\dot \gamma_{x,\xi}^i(t) \dot\gamma_{x,\xi}^j(t)  \phi(x)\, dt\, dx \\
  &= \int\int_{-\infty}^0  f_{ij}(y) \eta^i(y,\xi,t) \eta^j(y,\xi,t) \phi(y,\xi,t) J(y,\xi,t)\, dt\, dy,
\end{split}
\]
where $J$ is the corresponding Jacobian. Since $f$ is a distribution, our claim is proved. 

Therefore, \eqref{vv} is well defined in this case as well, and the rest of the proof remains the same.
\epr

\begin{remark} 
For $h$ as in \eqref{ha}, we therefore have $h_{in}=0$, $\forall i$,   everywhere in $ \tilde U$, and $h$ is supported in $M$ as well. 
Also, along all geodesics $\g$ lying in a small neighborhood of $x'=\mbox{const.}$, we have $If(\g) = Ih(\g)=0$. The so constructed $h$ depends on the choice of $\gamma_0$, and we will use it in the next section. Notice also that $h=0$ is an indication whether $f=dv$ in $\tilde U$ with some $v$ vanishing near one of the endpoints of $\gamma$ because the solution to \eqref{h2}, \eqref{h3} is unique. 
\end{remark}

We always work in $\tilde M$  below, therefore, in the next lemma, $\WF(f)$ is in the interior of $T^*\tilde M$.  We denote by $\pi$ the canonical projection onto the base. 

\begin{lemma}  \label{lemma_an1} Let $f\in \Ec'(\tilde M)$ be a symmetric tensor field in $\tilde M$ supported in $M$. 
Assume that $\WF(f)\cap \pi^{-1}(U%^{\rm{int}}
)$ does not contain covectors conormal to the geodesics $x'=\mbox{\rm const.}$, i.e, of the type $(x,\xi',0)$. Then the same is true for the tensor $h=f-dv$ constructed above. 
\end{lemma}

\bpr
It is enough to prove this for $v$. The proof for $v_n$ is straightforward since 
\begin{equation}\label{h1}
v_n(x) = \int_{-\infty}^{x^n} f_{nn}(x',y^n)\, dy^n.
\end{equation}
The integral above is a convolution with the Heaviside function $H$ and the analysis of its analytic wave front set follows easily, see e.g., sections~8.2 and 8.5 in \cite{H1}. We will not pursue this, instead we turn our attention to the proof for $v_i$, $i<n$, that implies the proof for $v_n$ as a partial case. 
Note that $\{v_i\}_{i=1}^{n-1}$, that we temporarily denote by $v$ (instead of $v'$) solves an ODE system of the type
\begin{equation}\label{h1a}
\partial_n v - A(x',x^n)v = w,\quad v|_{x^n\ll0 }=0,
\end{equation}
where   $A$ is a real-analytic matrix, $w$ satisfies the wave front assumptions of the lemma, and $w=0$ for $x^n<a(x')$. By the Duhamel's principle, the solution of \eqref{h1a} is given by
\[
v(x',x^n) = \int_{-\infty}^{x^n}\Phi(x',x^n,y^n)w(x',y^n)\, dy^n,
\]
where $\Phi$ is analytic. The r.h.s.\ above is an integral operator with kernel 
\[
K(x,y) = \Phi(x',x^n,y^n)H(x^n-y^n)\delta(x'-y')
\]
Its analytic singularities are conormal to the diagonal, i.e., they are included in the set
\[
\{(x,x,\xi,\eta); \; \xi+\eta=0\}.
\]
To analyze $\WF(v)$, we apply \cite[Theorem~8.5.5]{H1}. Note first that $\WF(K)_X$ is empty. In the notation of \cite[Theorem~8.5.5]{H1}, $\WF'(K)$ only changes the sign of $\eta$ and therefore, can contain covectors with $\xi=\eta$ only by the analysis above. Then $\WF'(K)\circ\WF(w)$, with $\WF'(K)$ considered as a relation, cannot contain covectors with $\xi_n=0$ because $w$ has the same property. 

Here $f$ is a distribution, as always. All arguments apply for distributions, as we showed in the proof of Lemma~\ref{lemma_v0}. Observe also that the kernel $K$ satisfies the conditions needed so that the corresponding integral operator extends to compactly supported distributions. 
\epr

\section{Analyticity along conormal directions}
The next proposition reflects the fact that $I$ and $\delta$ form a hypoelliptic system of operators.

\begin{prop}\label{S4:P1}
Let $\g_0$ be an open geodesic segment in $\tilde M$ with endpoints in $\tilde M \setminus M$. Let $f\in \Ec'(\tilde M)$ be a symmetric 2-tensor supported in $M$. 
Given  $(x_0,\xi^0)\in N^{*}\g_0\setminus 0$, assume that $(x_0,\xi^0)\not\in \WF(\delta f)$, and that $If(\g)=0$ for all $\g$ close enough to $\gamma_0$. Then $(x_0,\xi^0)\not\in \WF(f)$. 
\end{prop}

\bpr
We remark first that it is enough to assume that $If(\gamma)$ is real analytic, if $\gamma$ is parametrized in an analytic way. 

In the $C^\infty$ category, the proof follows directly from the following facts, see e.g., \cite{SU-AJM}. If $\psi$ is a standard cutoff that restricts the geodesics to a small neighborhood of $\gamma_0$, then the operator $I^*\psi I$ is a pseudo-differential operator elliptic on $N^*\gamma_0$ when applied to solenoidal tensors. In other words, if $Q$ is any elliptic pseudo-differential operator of order $-2$, then $(I^*I,Q\delta)$ forms a (matrix-valued) pseudo-differential operator of order $-2$ elliptic on $N^*\gamma_0$. Since $I^*\psi If=0$, and $\delta f$ is smooth at $(x,\xi)$, this completes the proof in the smooth case. 

The analytic case is more delicate because of the restrictions we have with the cut-offs. The pseudo-differential operator theory is well developed, see e.g., \cite{Treves} but its FIO analog is not so clear (see also \cite{Sj-A}). We follow here the approach based on the application of the analytic stationary phase method, see \cite{Sj-A}.

If $\delta f=0$, then the proposition still follows from \cite[Proposition~2]{SU-AJM}. When $\delta f$ is only microlocally analytic at $(x_0,\xi^0)$, it follows from the proof of that proposition, with one slight modification. Namely, we have to show that we still have, see equation (56) in \cite{SU-AJM}, that
\begin{equation}   \label{sol}
\int e^{\lambda [i(x-y)\cdot\xi-(x-y)^2/2]} C^j(x,y,\xi,\lambda)f_{ij}(x)\, dx = O(e^{-\lambda/C}), \quad \mbox{as $\lambda\to\infty$},
\end{equation}
with a classical analytic symbol $C^j$, in the sense of \cite{Sj-A}, having the following properties: $C^j$is  analytic near $(x_0,x_0,\xi^0)$ equal to $(\xi^0)^j$ when $x=y=x_0$, $\xi=\xi^0$; and equal to $0$ for $x$ outside some neighborhood of $x_0$. In order to do this, we start with
\[
\int e^{\lambda [i(x-y)\cdot\xi-(x-y)^2/2]} A(x,y,\xi)[\delta f]_i(x)\, dx = O(e^{-\lambda/C}), \quad \mbox{as $\lambda\to\infty$, $\forall i$},
\]
where $A$ is an appropriate cut-off elliptic near $(x_0,x_0,\xi^0)$, see \cite{Sj-A}. To get \eqref{sol}, we just integrate by parts. Then the proof of \cite[Proposition~2]{SU-AJM} works in our case as well without further modifications. 
\epr

\section{Proof of  Theorem~\ref{S1:T1}}\label{S4}
The ``if'' part of Theorem~\ref{S1:T1} is trivial. To prove the ``only if'' assertion of Theorem~\ref{S1:T1}(a), let  $\gamma_0\in\Ac$ be as in the theorem. We have that $\gamma_0$ can be continuously deformed within the set $\Ac$ to a point. Without loss of generality, we may assume that the only geodesic in that deformation, that is a point, is the final one. 
We extend all geodesics in $\Ac$ to maximal geodesics in  $\tilde M$ as before. Therefore, there exist two continuous curves $a(t)$, $b(t)$, $t\in[0,1]$, on $\PD \tilde M$ so that $\gamma_{[a(0),b(0)]}$ is tangent to $\PD M$ (and it is the only geodesic in the family $\gamma_{[a(t),b(t)]}$ with that property);  $\gamma_{[a(t),b(t)]}\in \Ac$ for any $t\in [0,1]$; and $\gamma_{[a(1),b(1)]}=\gamma_0$. It is enough to prove that $f=dv$ on $\tilde M_\Uc$ with $v=0$ outside $M$, where $\Uc$ is some neighborhood of $\gamma_0$.

Since $f$ and $f^s$ have the same line integrals in $M$, we can assume that $f=f^s$ in $M^{\inter}$, $f=0$ outside $M$; therefore $\delta f=0$ in $M^{\inter}$. For the argument that follows, we need to know that $f^s|_{\PD M}$ is analytic. To this end, we choose another analytic simple manifold $M_{1/2}$ so that $M\Subset M_{1/2}\Subset \tilde M$. Then $f$  is compactly supported in the interior of $M_{1/2}$. Let $f^s_{M_{1/2}}$ be the solenoidal projection of $f$ in $M_{1/2}$. Recall \cite{Sh1, SU2}, that $f^s_{M_{1/2}}=f-dv$, where $v$ solves the elliptic system $\delta dv = \delta f$ in ${M_{1/2}}$ with the regular boundary conditions $v=0$ on $\PD M_{1/2}$. Since $\delta f=0$ near $\PD M_{1/2}$, and $f$ is a distribution of finite order, the solution exists and is a distribution in some Sobolev space. Moreover, since  $\delta f=0$ near $\PD M_{1/2}$, we have that $v$ is analytic near $\PD M_{1/2}$, up to $\PD M_{1/2}$, see \cite{MorreyN}. Then so is $dv$. In $M_{1/2}\setminus M$, $f^s_{M_{1/2}}=-dv$, therefore, $f^s_{M_{1/2}}|_{\PD M_{1/2}}  =-dv|_{\PD M_{1/2}}$ is analytic.  
If we prove the theorem for $f^s_{M_{1/2}}$ in $M_{1/2}$, then we have the same for $f$ in $M$. Indeed, it is easy to show that $f=dv$ in ($M_{1/2})_\Ac$ with $v=0$ outside  $M_{1/2}$ implies $v=0$ in $M_{1/2}\setminus M$ by integrating \eqref{dv}. 
Now we can replace $M$ by $M_{1/2}$ in the proof that follows, and call in $M$ again. We also denote $f^s_{M_{1/2}}$ by $f$.

The advantage that we have now is that $f=f^s$ is solenoidal in $M$, and analytic near $\PD M$, up to $\PD M$. Let us denote by $v_0$ the restriction of $v$ to a collar neighborhood of $\PD M$, where $v_0$ is analytic. Therefore, after extending $M$ to $M_{1/2}$, and calling it $M$ again, we have the following.

\begin{lemma}  \label{lemma_v1}
There exists a neighborhood $V$  of $\PD M$ in $M$ of the kind $\mbox{\rm dist}(x,\PD M)<\epsilon_0$, $\epsilon_0>0$, 
and a uniquely defined $v_0$ there so that $f=dv_0$ in $V$, $v_0=0$ on $\PD M$, and $v_0$ is analytic in $V$, up to the boundary $\PD M$. 
\end{lemma}

\begin{proof}
Since we already proved existence of such $v_0$ near $\PD M$, we only need to show that  it is unique. Near any $p\in \PD_\Ac M$, introduce normal boundary coordinates $(x',x^n)$ as in the beginning of section~\ref{S2} but now $\partial/\partial x^n$ is normal to $\PD M$, and $x^n$ measures the distance to $\PD M$. Then we define $v_0$ as the solution of \eqref{h2}, \eqref{h3} for $x^n\ll1$. If we have that $f=dv$ near $p$ with a possibly different $v$ vanishing on $\PD M$,  we immediately get that $v=v_0$ because \eqref{h2}, \eqref{h3} has unique solution, and if $f=dv$, that solution is $v$. Now, $v_0$ is defined near $\PD M$, and $f=dv_0$. 
\end{proof}

We extend  $v_0$ as zero outside $M$. Then we still have $f=dv_0$ in $\tilde V := (\tilde M^{\inter}\setminus M) \cup V$. 
Let $p\in \gamma_{[a(0),b(0)]}$ be the point where  $\gamma_{[a(0),b(0)]}$ is tangent to $\PD M$. By Lemma~\ref{lemma_v1}, $f=dv_0$ near $p$. 

So  we get that the statement of the theorem is true in some neighborhood of the geodesics $\gamma_{[a(t),b(t)]}$, $0\le t\le 2t_0$, with some $0<t_0\ll1$ chosen so that all those geodesics are in $\tilde V$.

We want to propagate this property to a neighborhood of the  ``surface'' (that might be self-intersecting and non-smooth)
\[
\bigcup_{t\in [0,1]} \gamma_{[a(t),b(t)]}.
\]
For any $t$, define a small cone $\Cc_\epsilon(t)$ with vertex at $a(t)\in \PD \tilde M$  as follows. In the tangent space $T_{a(t)}\tilde M$, consider the cone of all vectors making angle less than $\epsilon>0$ with $\dot\gamma_{[a(t),b(t)]}$, where $\epsilon>0$ is fixed small enough.  Then we define $\Cc_\epsilon(t)$
as the image of that cone in $\tilde M$ under the exponential map. The choice of $\epsilon$ is the following: we require that 
\begin{itemize}
 \item[(i)] $\Cc_{2\epsilon}(t)\subset\tilde  M_{\Ac}, \quad \forall t\in [0,1]$;
 \item[(ii)] none of the geodesics inside the cone $\bar\Cc_{2\epsilon}(t)$, $t_0\le t\le1$ with vertices at $a(t)$ is tangent to $\PD M$;  
 \item[(iii)]  $\Cc_\epsilon(t)\subset \tilde V$ for $0,\le t\le t_0$. 
\end{itemize}
This   can be arranged by a compactness argument. 

For any $t$, construct a tensor field $h_t$ in $C_{2\epsilon}(t)$ as in the remark following  the proof of Lemma~\ref{lemma_v0}. Then $h_t=0$ outside $M$ by Lemma~\ref{lemma_v0}. For $t\le t_0$ we have that $h_t=0$ in  $\Cc_{\epsilon}(t)$ by (iii).  
Set
\[
t^* = \sup \left\{  t\in (0,1]:\; \mbox{$h_t=0$ in $\Cc_\epsilon(t)$}\right\}.
\]
We aim to show that $t^*=1$. Assume that $t^*<1$. Then one can show that  $h_{t^*}=0$ in $\Cc_\epsilon(t^*)$ because $h_{t^*}=0$ outside $M$. 

We will show now that $h_{t^*}=0$ in  $\Cc_{2\epsilon}(t^*)$ as well. That will yield a contradiction because then one can increase $t^*$ slightly to $t$ to get $\Cc_\epsilon(t)\cap M \subset \Cc_{2\epsilon}(t^*)\cap M$ for all $t>t^*$ close enough to $t^*$.  That would contradict the choice of $t^*$. 

To fulfill this program,  consider $h_{t^*}$ in  $\Cc_{2\epsilon}(t^*)$.  The support of $ h_{t^*}$ is included in $M$ and it does not intersect the interior of $\Cc_\epsilon(t^*)$. 
Increase the angle from $\epsilon$ to $2\epsilon$ until it does (if it does not, we are done). Let $\epsilon_0$ with $\epsilon<\epsilon_0\le2\epsilon$ be the smallest number with that property. 

Consider the cone $\Cc_{\epsilon_0}(t^*)$. Then $h_{t^*}=0$ in $\Cc_{\epsilon_0}(t^*)$, and $\mbox{supp}\, h_{t^*}$ and $\Cc_{\epsilon_0}(t^*)$ have a common point $q$ lying on the boundary of each set. The point $q$ cannot be on $\PD \tilde M$ because $h_{t^*}=0$ outside $M$. 
 So $q$ is an interior point of $\tilde M$.

We have that in $\tilde M$, $(\delta f)_i= f^s_{ij}\nabla^j\chi = -f^s_{ij}\nu^j\delta_{\PD M}$, where $\chi$ is the characteristic function of $M$. Therefore, $\delta f$ may have only conormal analytic singularities at $\PD M$, i.e., singularities in $N^*\PD M$. 
Let $\gamma$ be the geodesic in $\tilde M$ on the surface $\partial\Cc_{\epsilon_0}(t^*)$ that contains $q$. 
Since $N^*\gamma$ does not intersect $N^*\PD M$ by (ii), we get that $f$ has no analytic singularities in $N^*\gamma$ by Proposition~\ref{S4:P1}. We want to point out that this is true even if the base point is on $\PD M$. By Lemma~\ref{lemma_an1}, $h$ has no analytic singularities in $N^*\gamma$ either. 

Therefore, $h_{t^*}=0$ in $\Cc_{\epsilon_0}(t^*)$, and $h_{t^*}$ is defined at least in a small enough neighborhood of $q\in \partial \Cc_{\epsilon_0}(t^*)$. Since the conormals to $\Cc_{\epsilon_0}(t^*)$ at $q$ are not in $\WF(h)$, we have that $h_{t^*}=0$ near $q$ by the Sato-Kawai-Kashiwara theorem, see e.g., \cite{SKK} or \cite[Theorem~8.2]{Sj-A}. This contradicts the fact that $q$ is on the boundary of $\mbox{supp}\,h_{t^*}$. Therefore, $t^*=1$. 
This completed the proof of Theorem~\ref{S1:T1}(a).

Consider part (b): In any open set, where $f^s=dv$, we have $\delta d v=0$. Therefore, $v$ is analytic there, by elliptic regularity. Thus in any cone $\Cc_\epsilon(t)$ as above, the field $v$ having the property that $f=dv$ in $\Cc_\epsilon(t)$, is in fact an analytic continuation of $v_0$ from a neighborhood of $\PD M$ in $M$ to $\Cc_\epsilon(t)\cap M$. 
If  $M_{\Ac}$ is simply connected, then we have  uniqueness of that continuation.
Under the assumption \eqref{ii}, this still works, as we show below, see also \cite{SU-AJM} for a similar argument. 

Let $q\in M_\Ac$, and let $p\in \PD_\Ac M$ be as in the proof above. Having a path connecting $q$ and $p$, we define $v$ near $q$ as an analytic continuation of $v_0$ from the neighborhood $V$ of $p$, see Lemma~\ref{lemma_v1},  along that path. We need to show that this definition is independent of the choice of the path. Any other path is homotopic to the composition of the first one and a path on $\PD_\Ac M$, by \eqref{ii}. Near the boundary path, $f=dv_0$ by Lemma~\ref{lemma_v1}. This allows us to show that the two analytic continuations coincide.

This completes the proof of Theorem~\ref{S1:T1}.

\section{Proof of  Theorem \ref{S1:T2}}

We now prove Theorem \ref{S1:T2}. The proof of this theorem follows from the following lemmas and Theorem \ref{S1:T1}.

\begin{lemma}\label{S2:L1}
For any  $x\in M\setminus K$, there is a geodesic passing through $x$ that does not intersect $K$.
\end{lemma}
\bpr
Consider $\xi\in S_xM$. If the geodesic starting at $x$ in the direction $\xi$ intersects $K$, then the other half of this geodesic segment, that is, the one starting at $x$ in the direction $-\xi$ does not intersect $K$ by the geodesic convexity of $K$. So if $L$ denotes the set of unit directions $\xi$ such that the geodesic passing through $x$ in the direction $\xi$ intersects $K$, then we have $(-L)\cap L=\emptyset$, where $-L$ is the set of unit vectors $\xi$ such that $-\xi\in L$. If we assume that every geodesic passing through $x$ intersects $K$, then $(-L)\cup L=S_{x}M$. The set $L$ is closed because, if a geodesic say $\g$ passing through $x$ in the direction $\xi \in S_{x}M$ does not intersect $K$, there is a neighborhood of $\xi$ in $S_{x}M$ such that the geodesics in these directions also do not intersect $K$. We now have a disconnection of $S_{x}M$  which is a contradiction.
\epr
\begin{lemma}\label{S2:L2}
Let $x\in M\setminus K$. Then there is a continuous deformation of a geodesic passing through $x$ that does not intersect $K$ by such  geodesics  to a point on the boundary $\PD M$.
\end{lemma}
\bpr
By Lemma~\ref{S2:L2}, there is a geodesic $\gamma_0$ through $x$ lying in $\tilde M\setminus K$. There is $\epsilon>0$ so that the cone $\Cc_0$ constructed as in the previous section, with vertex at one of the endpoints of $\gamma_0$ and angle $\epsilon$, still does not intersect $K$. 
Consider any continuous deformation of that cone to a cone in $\tilde M\setminus M$. If none of the cones in this deformation touches $K$, then we are done. If some of them does, let us denote the first cone in that family by $\Cc_1$. Let $\gamma_1$ be its axis. It is enough to show that $\gamma_1$ can be deformed to a geodesic in $\tilde M\setminus M$ by a homotopy that does not intersect $K$.

Assume that $n\ge3$ first. 
Let $p\in \Cc_1\cap K$. Then $p\in \PD \Cc_1$. Let $S$ be the collection of all maximal geodesics in $\tilde M$ through $p$ tangent to $C_1$. Then $S$ is a smooth surface that is the image of $T_p\Cc_1$ under the exponential map; more precisely, $S = \exp_p \big( \exp_p^{-1}(\tilde M)\cap  T_p\Cc_1   \big)$. Then $S$ separates $\tilde M$ into two parts $\tilde M_1$ and $\tilde M_2$ (we do not include $S$ in either of them, so actually, $\tilde M=\tilde M_1\cup S\cup \tilde M_2$). This is easily seen by considering $\exp^{-1}(\tilde M)$ first. In particular,  $S$ is diffeomorphic to an $(n-1)$-dimensional ``disk'', i.e., to an $(n-1)$-dimensional ball. On the other hand, $\tilde M$ is diffeomorphic to an $n$-dimensional ball.  

We will show that $K$ is included in the closure of one of the $\tilde M_i$, and does not intersect the other one. Assume that this is not true. Then there are two points $q_1$ and $q_2$ in $K$, so that $q_i$ belongs to the interior of $M_i$, $i=1,2$. Then the geodesics $\gamma_{[p,q_i]}$, $i=1,2$, are contained in $K$ by the convexity assumption.  
Moreover, their only common point with $S$ is the endpoint $p$, by the simplicity condition. Therefore, $\gamma_{[p,q_i]} \in (M_i\cup S)\cap K$,  $i=1,2$. 
Then the velocity  vectors  of $\gamma_{[p,q_i]}$ at $p$ (with the orientation fixed to be from $p$ to $q_i$) cannot be tangent to $T_pS$ and belong to different half-spaces of $T_p\tilde M$. Since $T_p\PD \Cc_1= T_pS$,
this shows that a small enough part of one of those geodesics near $p$ is inside the cone $\Cc_1$, and a small part of the other is outside it. This contradicts the fact that the interior of $\Cc_1$ does not contain points  of $K$.  

We will deform now $\gamma_1$ to a geodesic in $\tilde M\setminus M$, without intersecting $K$. Let $\tilde M_1$ be the set  which contains $K$. Then $\gamma_1\in \tilde M_2$. The part of the boundary $\PD \tilde M \cap \tilde M_2$ is diffeomorphic to a ``disk'', i.e., to an $(n-1)$-dimensional ball. In particular, it is connected. Since $\gamma_1\in \tilde M_2$, we fix one of the endpoints of $\gamma_1$ and move the other along a smooth curve in  $\PD \tilde M \cap \tilde M_2$  until it reaches the first one. 

If $n=2$, then $\Cc$ consists of two geodesics, and $\gamma_1$ is one of them. Then the proof above works in the same way with $\gamma_1$ playing the role of $S$. 
\epr

\begin{lemma}\label{S4:L3}
$\pi_{1}(M\setminus K, \PD M)=0$.
\end{lemma}
\bpr
Let  $c$ be a path in $M\setminus K$ with a base point on $\PD M$.
Fix a point $p\in K$ and let $x\mapsto \mbox{proj}_p(x)\in \PD M$ be the projection that maps $x$ to the endpoint of the geodesic from $p$ to $x$ until it hits $\PD M$.  Using $\mbox{proj}_p(x)$, project $c$ on $\PD M\cong \Sb\sp{n-1}$. This provides a way to continuously deform $c$ to its projection on $\PD M$. 
\epr

\noindent {\it Proof of Theorem \ref{S1:T2}}\,: Now the lemmas above together with Theorem \ref{S1:T1} prove Theorem \ref{S1:T2}.

\bigskip
\textbf{Acknowledgments.} 
The first author would like to thank John Lee for the several discussions on this work. Both authors thank Gunther Uhlmann for suggesting this problem.

%================

%=============

\end{document}